\documentclass{amsart}
\usepackage[T1]{fontenc}
\usepackage{amsmath,amssymb,amsthm}
\usepackage{mathrsfs}
\usepackage[a4paper,margin=2.5cm]{geometry}
\usepackage[colorlinks,linkcolor=blue,citecolor=blue,urlcolor=blue]{hyperref}
\usepackage{amssymb}
\usepackage{dsfont}   % 必须加载这个包
\newtheorem{theorem}{Theorem}[section]
\newtheorem{lemma}[theorem]{Lemma}
\newtheorem{proposition}[theorem]{Proposition}

\title[Recurrence and Uniqueness for Markov Processes]{Recurrence and Uniqueness Criteria for Markov Processes from Moments of Integral-Type Functionals}

\author{Ling-Di Wang}
\address{School of Mathematics and Statistics, Henan University, Kaifeng 475001, China}
\email{wanglingdi@henu.edu.cn}
\author{Yu-Hui Zhang}
\address{School of Mathematical Sciences,  Beijing Normal University, Beijing,100875, China}
\email{zhangyh@bnu.edu.cn}
\author{Zu-Yi Zhao}
\address{School of Mathematical Sciences, Beijing Normal University, Beijing, 100875, China}
\email{zhaozy@mail.bnu.edu.cn}
\thanks{Corresponding author: zhaozy@mail.bnu.edu.cn(Zu-Yi Zhao)}

\begin{document}

\begin{abstract}

We investigate recurrence and uniqueness of continuous-time Markov chains using moments of integral-type functionals. Unlike the classical criteria from Chen (2004) and the Feller–Reuter criterion, we develop general unified conditions for recurrence and uniqueness of Markov processes based on moments of integral-type functionals. The associated criteria follow directly from these conditions whenever the mean functionals can be expressed explicitly. As applications, explicit recurrence and uniqueness criteria are derived for single-birth, single-death, and birth-death processes on trees.\end{abstract}
\maketitle

\vspace{0.5em}

\section{Background and Main Results}
The uniqueness and recurrence of Markov processes determined by a $Q$-matrix are fundamental problems in the theory of continuous-time Markov chains. Let $(q_{ij}:i,j\in E)$ be a totally stable and conservative transition rate matrix (also called a $Q$-matrix) on a countable state space $E:=\{0,1,2,\dots\}$. That is, $q_{ij}\ge 0$ for $j\neq i$, $q_i:=-q_{ii}\ge 0$, $q_i<+\infty$ for all $i\in E$, and $q_i=\sum_{j\neq i} q_{ij}$ for $i\in E$.  The $Q$-matrix $Q$ is called irreducible if for any pair of states $i,j\in E$, there exists a finite sequence of states $i=k_0,k_1,\dots,k_m=j$ such that $q_{k_{l},k_{l+1}}>0$ for each $0\le l\le m-1$. Every $Q$-matrix admits a minimal continuous-time Markov process $\{X^{\text{min}}(t):t\ge 0\}$, with transition probabilities $P^{\text{min}}(t)=(P_{ij}^{\text{min}}(t):i,j\in E)$. 
We say $Q$-matrix is  regular if $\{X^{\text{min}}(t)\}$ is the only $Q$-process associated with $Q$.
 For conservative $Q$, irreducibility of $Q$ implies irreducibility of the minimal process: for any $i,j\in E$, $P_{ij}^{\text{min}}(t)>0$ for some $t>0$.
An irreducible minimal process is recurrent if
\[
\int_0^\infty P_{00}^{\text{min}}(t)\,\mathrm dt=+\infty,
\]
and transient otherwise.

Most available criteria for recurrence and uniqueness rely on analytical constructions or drift arguments. In contrast, our method uses moments of integral-type functionals. We first recall two classical criteria for $Q$-processes below, which serve as a useful reference.

Building on Feller’s work on the recurrence of minimal processes (\cite{F-1957}), Chen derived the recurrence criterion given below.The uniqueness criterion we recall is the classical Feller-Reuter criterion (\cite{F-1957,R-1957}). Let $Q=(q_{ij})$ be a totally stable, conservative and irreducible $Q$-matrix.
\begin{itemize}
\item The minimal process $P^{\text{min}}(t)$ is recurrent if and only if the system
\[
x_i = \sum_{j\ne i} \frac{q_{ij}}{q_i}x_j,\quad 0\le x_i \le 1,\ i\in E
\]
admits only the trivial solution; see \cite[Lemma 4.51]{Chen-2004}.

\item The $Q$-process determined by $Q$ is unique if and only if for some (equivalently, for all) $\lambda>0$, the system
\[
u_i=\sum_{k\ne i}\frac{q_{ik}}{\lambda+q_i}u_k,\quad 0\le u_i\le 1,\ i\in E
\]
admits only the trivial solution $u_i\equiv 0$ (\cite{F-1957,R-1957}).
\end{itemize}
Explicit recurrence criteria for birth-death and single-birth processes can follow directly from \cite[Lemma 4.51]{Chen-2004}. By contrast, deriving such explicit criteria for single-death processes remains challenging.
Criteria for uniqueness of birth-death and single-birth processes can follow from the Feller-Reuter theorem, though the derivation is nontrivial. For single-death processes, the task of extracting explicit uniqueness criteria becomes considerably more difficult(refer to \cite{M-Y-Z-2026}). Different from these classical criteria, we adopt a distinct discriminative strategy: we derive general unified conditions for recurrence and uniqueness of Markov processes in terms of moments of integral functionals, from which the criterion is fully characterized by the limits of moments of the hitting time.

Before stating our main results, we introduce some notation and preliminaries. For a nonempty subset $A\subset E$, define the return time to the set $A$ of the process:
\begin{equation*}
\sigma_A = \inf\{t \ge \eta_1: X^{\min}(t)\in A\}.
\end{equation*}
where $\eta_1$ denotes the first jump time of $\{X^{\min}(t)\}$. Let $V$ be a nonnegative and nonvanishing function on $E$. The integral-type functional is defined as \(\int_{0}^{\sigma_A} V(X^{\min}(t))dt\). We denote the moment of  integral-type functional as
\[
T_{i,A}(V) := \mathbb{E}_i\int_{0}^{\sigma_A} V(X^{\min}(t))dt,
\]
It is known from  \cite{H-G-1988} that
%\begin{theorem}(\cite{H-G-1988})
$(T_{i,A}(V),\,i\in E)$ is the minimal nonnegative solution to the following system of equations:
\[
x_i = \sum_{j\notin A\cup\{i\}} \frac{q_{ij}}{q_i} x_j + \frac{V(i)}{q_i},\quad i\in E.
\]
When $A=\{n,n+1,n+2,\dots\}$ for some positive integer $n$, we write $T_{i,A}(V)$ and $\sigma_A$ as $T_{i,n+}(V)$ and $\sigma_{n+}$, respectively. When $V(k)=\mathds 1_{\{k=0\}}=\delta_{0,k}$, we write $T_{i,n+}(V)$ as $T_{i,n+}$. By the localization theorem, $(T_{i,n+}(V),\,0\le i<n)$ is the minimal nonnegative solution of
\begin{equation}\label{eq0.1}
x_i = \sum_{\substack{j\ne i\\j<n}} \frac{q_{ij}}{q_i} x_j + \frac{V(i)}{q_i},\quad 0\le i < n,
\end{equation}
Denote by $I^{(n)}$ the $n\times n$ identity matrix and by $\Pi^{(n)}=(\Pi_{ij}^{(n)})_{0\le i,j<n}$ the associated locally embedded chain with
\[
\Pi_{ij}^{(n)}:= \frac{q_{ij}}{q_i}\,\mathds{1}_{\{j\ne i\}}.
\]
Let
\[
\boldsymbol{x}^{(n)} = (x_0,x_1,\dots,x_{n-1})^\top,\qquad
\boldsymbol{r}^{(n)}=\left(\frac{V(0)}{q_0},\frac{V(1)}{q_1},\dots,\frac{V(n-1)}{q_{n-1}}\right)^\top
\]
be $n$-dimensional column vectors. Then \eqref{eq0.1} admits the compact matrix form
\[
\big({I}^{(n)} - {\Pi}^{(n)}\big)\boldsymbol{x}^{(n)} = \boldsymbol{r}^{(n)}.
\]

By irreducibility of $Q$, this finite-state locally embedded chain is  transient. Hence its Green matrix is entrywise finite:
\[
\big({I}^{(n)} -{ \Pi}^{(n)}\big)^{-1} = \sum_{m=0}^{\infty} \big({\Pi}^{(n)}\big)^m,
\]
where the series converges entrywise. It follows that system \eqref{eq0.1} admits a unique solution
\[
\boldsymbol{x}^{(n)} = \sum_{m=0}^{\infty} \big({\Pi}^{(n)}\big)^m \boldsymbol{r}^{(n)}.
\]
Specially, when $V(k)=\delta_{0,k}$, this solution may also be written via Cramer’s rule as
\begin{equation}\label{eq0.2}
T_{k,n+} = \frac{\big|{\Delta}_{k+1}^{(n)}\big|}{\big|{I}^{(n)}-{\Pi}^{(n)}\big|},\quad k=0,1,2,\dots,n-1,
\end{equation}
where ${\Delta}_{j}^{(n)}$ denotes the matrix obtained by replacing the $j$-th column of ${I}^{(n)}-{\Pi}^{(n)}$ with the column vector $r^{(n)}$, and $|{A}|$ stands for the determinant of matrix ${A}$.

In this paper, our main focus is the explicit characterization of $(T_{i,n+},\,0\le i<n)$, which have been explained above. As will be shown, these quantities yield necessary and sufficient conditions for the recurrence and uniqueness of Markov processes. For ease of presentation, we write $\operatorname{M-lim}_{n\to\infty} a_n$ to denote the (generalized) limit of a monotone sequence $\{a_n\}$ throughout the paper.

Our first main result establishes a recurrence criterion for the irreducible Q-process stated below:
\begin{theorem}[Recurrence Criterion]\label{Rec}
Let $Q=(q_{ij})$ be a totally stable, conservative and irreducible $Q$-matrix. Then $P^{\text{min}}(t)$ is recurrent if and only if
\[
\operatorname{M-lim}_{n\to\infty} T_{0,n+}=+\infty.
\]
%where $V(k)=\mathbf{1}_{\{0\}}(k)=\delta_{0k}$.
\end{theorem}

 Given a totally stable and conservative $Q$-matrix ${Q}=(q_{ij})$, we define a new $Q$-matrix on $\overline{E}:=E\cup\{-1\}$ by
\begin{equation}\label{eq:newQ}
\overline{{Q}}=
\begin{pmatrix}
-1 & \boldsymbol{p}\\
\boldsymbol{1}^\top & {Q}-{I}
\end{pmatrix}.
\end{equation}
where $\boldsymbol{p}=(p_0,p_1,p_2,\dots)$ with $p_i\ge 0$ and $\sum_{i\ge 0}p_i=1$. The values of $p_i$ are chosen merely to ensure the irreducibility of $\overline{Q}$.
Let $\boldsymbol{1}=(1,1,1,\dots)$ and ${I}=(\delta_{ij}:i,j\in E)$. Denote by $\{\overline{X}^{\text{min}}(t)\}$ (or $\overline{P}^{\text{min}}(t)$) the minimal process corresponding to $\overline{Q}$, and write $\bar{\sigma}_{n+}$ for the return time to the set $\{n,n+1,n+2,\dots\}$.
Define
\[
\overline{T}_{-1,n+}=\mathbb{E}_{-1}\int_{0}^{\bar{\sigma}_{n+}} \mathds{1}_{\{\overline{X}^{\text{min}}(t)=-1\}}\mathrm dt,
\]
which is just $T_{i,A}^{(1)}(V)$ with $V(k)=\delta_{-1,k}$ and $A=\{n,n+1,n+2,\dots\}$. 
Our second main result establishes a uniqueness criterion for the Q-process stated below.
\begin{theorem}[Uniqueness Criterion]\label{uni}
Let $Q=(q_{ij})$ be a totally stable and conservative $Q$-matrix. The $Q$-process determined by $Q$ is unique if and only if
\[
\operatorname{M-lim}_{n\to\infty}\overline{T}_{-1,n+}=+\infty.
\]
\end{theorem}

For a $Q$-matrix, once $T_{0,n+}$ or $\overline{T}_{-1,n+}$ is computed explicitly, Theorem \ref{Rec} or Theorem \ref{uni} yields explicit criteria for recurrence or uniqueness, respectively. We illustrate our results with examples of single-birth processes, single-death processes, and birth-death processes on trees. Our method applies to uniformly finite-range jump Markov processes, and our results are compared with those in \cite{Z-2023}. Furthermore, the approach extends to recurrence criteria for discrete-time Markov chains by replacing the integral functional with an additive functional. We omit further details for brevity.

The rest of the paper is organized as follows. Section 2 is devoted to the proofs of our main results. In Sections 3, 4 and 5, we apply our results to single-birth processes, single-death processes, and birth-death processes on trees, respectively.
\medskip

To conclude this section, we discuss an alternative probabilistic interpretation of $T_{0,n+}$.
Define the hitting time to a set $A$ as
\[
\tau_A = \inf\{t \ge 0: X^{\text{min}}(t)\in A\},
\]
and let $\tau_{n+}=\tau_A$ with $A=\{n,n+1,n+2,\dots\}$.
By the strong Markov property, $(\mathbb{P}_i(\tau_{n+}<\sigma_0),\,0\le i<n)$ is the unique solution to the linear system
\begin{equation}\label{eq0.3}
x_i=\sum_{\substack{0<j<n\\j\ne i}} \frac{q_{ij}}{q_i}\,x_j+\frac{\sum_{k\ge n} q_{ik}}{q_i},\qquad 0\le i<n.
\end{equation}
Let $I^{(n)} -\widetilde{\Pi}^{(n)}$ be the coefficient matrix of \eqref{eq0.3}.
 and let ${\boldsymbol{\widetilde{r}}^{(n)}}$ denote the $n$-dimensional column vector
\[
\Big(\frac{\sum_{k\ge n}q_{0k}}{q_0},\,\frac{\sum_{k\ge n}q_{1k}}{q_1},\,\dots,\frac{\sum_{k\ge n}q_{n-1,k}}{q_{n-1}}\Big)^T.
\]
By Cramer’s rule,
\[
\mathbb{P}_0(\tau_{n+}<\sigma_0)=\frac{|\widetilde{\Delta}_1^{(n)}|}{|I^{(n)} -\widetilde{\Pi}^{(n)}|},
\]
where $\widetilde{\Delta}_1^{(n)}$ is the matrix obtained by replacing the first column of $I^{(n)} -\widetilde{\Pi}^{(n)}$ with the column vector $\boldsymbol{\widetilde{r}^{(n)}}$.

Now take $V(k)=\delta_{0,k}$ in system \eqref{eq0.1}, by \eqref{eq0.2}, we get 
\[
T_{0,n+}=\frac{|\Delta_1^{(n)}|}{|I^{(n)} -{\Pi}^{(n)}|}.
\]
where $\Delta_1^{(n)}$ is the matrix obtained by replacing the first column of $I^{(n)}-{\Pi}^{(n)}$ with the vector $\big(\tfrac{1}{q_0},0,\dots,0\big)^T$.

It is easy to see that $|\Delta_1^{(n)}|=|I^{(n)} -\widetilde{\Pi}^{(n)}|/q_0$ and $|I^{(n)}-\Pi^{(n)}|=|\widetilde{\Delta}_1^{(n)}|$. Combining the preceding discussion, we arrive at
\[
T_{0,n+}=\frac{1}{q_0\mathbb{P}_0(\tau_{n+}<\sigma_0)}.
\]
Taking the limit as $n\to\infty$, we  obtain
\[
\int_{0}^{+\infty}P_{00}^{\text{min}}(t)\mathrm dt=\frac{1}{q_0\mathbb{P}_0(\eta\le\sigma_0)}.
\]
where $\eta=\lim_{n\to\infty}\tau_{n+}$ is the explosion time of the Q-process.

\section{Proofs of the Main Results}
{\noindent{ \bf Proof of Theorem \ref{Rec}}} Let $\eta=\text{M-}\lim_{n\to\infty}\sigma_{n+}=\lim_{n\to\infty}\tau_{n+}$ be the explosion time of $P^{\text{min}}(t)$. 
By the monotone convergence theorem,
\begin{align*}
\int_{0}^{+\infty} P_{00}^{\text{min}}(t)\mathrm dt
&= \int_{0}^{+\infty} \mathbb P_0(X^{\text{min}}(t)=0)\mathrm dt \\
&= \int_{0}^{+\infty} \mathbb E_0 \mathds 1_{\{X^{\text{min}}(t)=0,\, t\le \eta\}} \mathrm dt \\
&= \lim_{n\to\infty}\mathbb E_0 \int_{0}^{\sigma_{n+}} \mathds 1_{\{X^{\text{min}}(t)=0\}} \mathrm dt \quad(\text{monotone increasing limit})\\
&= \text{M-}\lim_{n\to\infty} T_{0,n+}.
\end{align*}
 Combining this with the definition of recurrence, the conclusion of Theorem \ref{Rec} follows. $\square$

\vskip0.2cm 
Before proving Theorem \ref{uni}, we state the following lemma, originally due to Feller (1957) .
\begin{lemma}(\cite{F-1957})\label{Feller}
Let $Q=(q_{ij})$ be a totally stable and conservative $Q$-matrix. Then
\[
\int_{0}^{\infty} P_{ij}^{\text{min}}(t)\mathrm dt
= \frac{1}{q_j}\sum_{m=0}^{\infty} \big(\Pi^m\big)_{ij},\qquad
\Pi_{ij}= \frac{q_{ij}}{q_i}\mathds{1}_{\{i\ne j\}}.
\]
Suppose further that $Q$ is irreducible. Then the minimal process $\{P^{\text{min}}(t)\}$ is recurrent if and only if its embedded chain is recurrent.
\end{lemma}

 References \cite{Chen-2015,S-2015} establish the relation between the uniqueness of the original process and the recurrence of the minimal process on the extended state space. Based on Lemma \ref{Feller}, we reformulate their result as the following lemma.

\begin{lemma}\label{rec-lem-1}
Let $Q=(q_{ij})$ be a totally stable and conservative $Q$-matrix. Define $\overline{Q}$ on $\overline{E}:=E\cup\{-1\}$ by \eqref{eq:newQ}. Then the process determined by $Q$ is unique if and only if the minimal process corresponding to $\overline{Q}$ is recurrent.
\end{lemma}

\noindent{\bf Proof of Theorem \ref{uni}}
Take $V(k)=\delta_{-1,k}$. Let $\overline{X}^{\text{min}}(t)$ denote the minimal process corresponding to $\overline{Q}$, with $\overline{\sigma}_{n+}$ its return time and $\overline{P}^{\text{min}}(t)$ its transition matrix. From the argument in the proof of Theorem \ref{Rec}, we obtain
\[
\int_{0}^{+\infty}\overline{P}_{-1,-1}^{\text{min}}(t)\mathrm dt = \text{M-}\lim_{n\to\infty}\overline{T}_{-1,n+}.
\]
Hence the desired conclusion follows immediately from Lemma \ref{rec-lem-1} and Theorem \ref{Rec}. $\square$

\section{ Single Birth Processes}
A Markov process is called a single birth process if $q_{i,i+1}>0$ and $q_{i,i+j}=0$ for all $i\in E$, $2\le j\in E$. Define $q_n^{(k)}=\sum_{j=0}^{k} q_{nj}$ for $0\le k<n$, and set
\[
F_i^{(i)}=1,\qquad F_n^{(i)}=\frac{1}{q_{n,n+1}}\sum_{k=i}^{n-1} q_n^{(k)} F_k^{(i)},\quad n>i.
\]
For the subsequent proofs, we state the following lemma from \cite{C-Z-2014}.
\begin{lemma}(\cite[Corollary 2.3 ]{C-Z-2014})\label{C-Z-1}
For a given function $f$, the sequence $\{h_n\}$ defined recursively by
\[
h_n=\frac{1}{q_{n,n+1}}\bigg(\sum_{i\le k\le n-1} q_n^{(k)} h_k + f_n\bigg),\quad n\ge i,
\]
has the unified representation
\[
h_n=\sum_{k=i}^{n} \frac{F_n^{(k)} f_k}{q_{k,k+1}},\quad n\ge i.
\]
Moreover, $\{F_n^{(k)}\}$ satisfies
\[
F_i^{(i)}=1,\qquad F_n^{(i)}=\sum_{k=i+1}^{n} \frac{F_n^{(k)} q_k^{(i)}}{q_{k,k+1}},\quad n>i.
\]
\end{lemma}

We have the following explicit result for the system \eqref{eq0.1} associated with single birth processes.

\begin{proposition}\label{B-f1}
For a single birth $Q$-matrix, system \eqref{eq0.1} has a unique solution
\[
x_i=\sum_{i\le \ell\le n-1}\sum_{0\le k\le \ell} \frac{F_\ell^{(k)} V(k)}{q_{k,k+1}},\quad 0\le i\le n-1.
\]
\end{proposition}

\begin{proof}
By \eqref{eq0.1} and the conservativeness of $Q$, for $0\le i\le n-1$, we obtain
\[
\begin{aligned}
\mathds 1_{\{i=n-1\}} q_{i,i+1}x_i + \mathds 1_{\{0\le i\le n-2\}} q_{i,i+1}(x_i-x_{i+1})
&=\sum_{0\le j\le i-1} q_{ij}(x_j-x_i)+V(i) \\
%&=\sum_{0\le j\le i-1} q_{ij}\sum_{j\le \ell\le i-1}(x_\ell-x_{\ell+1})+V(i) \\
&=\sum_{0\le \ell\le i-1} q_i^{(\ell)}(x_\ell-x_{\ell+1})+V(i).
\end{aligned}
\]
Therefore, for all $0\le i\le n-2$, it follows from Lemma \ref{C-Z-1} that
\begin{equation*}\label{Birth-f1}
x_i - x_{i+1} = \sum_{0\le k\le i} \frac{F_i^{(k)}V(k)}{q_{k,k+1}}.
\end{equation*}
For $i=n-1$,
\[
\begin{aligned}
q_{n, n-1}x_{n-1} &= \sum_{0\le \ell\le n-2} q_{n-1}^{(\ell)}(x_\ell - x_{\ell+1}) + V(n-1) \\
&= \sum_{0\le \ell\le n-2} q_{n-1}^{(\ell)}\sum_{0\le k\le \ell} \frac{F_\ell^{(k)}V(k)}{q_{k,k+1}} + V(n-1) \\
&= \sum_{0\le k\le n-2} \frac{V(k)}{q_{k,k+1}} \sum_{k\le \ell\le n-2} q_{n-1}^{(\ell)}F_\ell^{(k)} + V(n-1) \\
&= q_{n-1,n}\sum_{0\le k\le n-2} \frac{F_{n-1}^{(k)}V(k)}{q_{k,k+1}} + V(n-1) .\\
%&= q_{n-1,n}\sum_{0\le k\le n-1} \frac{F_{n-1}^{(k)}V(k)}{q_{k,k+1}}.
\end{aligned}
\]
Hence, 
\[
x_{n-1} = \sum_{0\le k\le n-1} \frac{F_{n-1}^{(k)}V(k)}{q_{k,k+1}}.
\]
Therefore, for $0\le i\le n-1$,
\[
\begin{aligned}
x_i &= \sum_{i\le \ell\le n-2} (x_\ell - x_{\ell+1}) + x_{n-1} \\
&= \sum_{i\le \ell\le n-2}\sum_{0\le k\le \ell} \frac{F_\ell^{(k)}V(k)}{q_{k,k+1}} + \sum_{0\le k\le n-1} \frac{F_{n-1}^{(k)}V(k)}{q_{k,k+1}} \\
&= \sum_{i\le \ell\le n-1}\sum_{0\le k\le \ell} \frac{F_\ell^{(k)}V(k)}{q_{k,k+1}}.
\end{aligned}
\]
The proof  is complete. \end{proof}

Fix some $j\in \{0, 1, 2,\cdots, n-1\}$. Take $V(k)=\mathds 1_{\{k=j\}}$. Then by Proposition \ref{B-f1},
\[
x_i=
\begin{cases}
\displaystyle\sum_{i\le \ell\le n-1}\frac{F_\ell^{(j)}}{q_{j,j+1}}, & 0\le j\le i\le n-1;\\[8pt]
\displaystyle\sum_{j\le \ell\le n-1}\frac{F_\ell^{(j)}}{q_{j,j+1}}, & 0\le i<j.
\end{cases}
\]
In particular, for $j=0$, we have $T_{0,n+}=x_0=\sum_{0\le \ell\le n-1}F_\ell^{(0)}/q_{01}$. Combining this with Theorem \ref{Rec}, we immediately obtain the following criterion (cf. \cite{Chen-1999, Y-C-1986}).
\begin{proposition}
Let $Q=(q_{ij})$ be a totally stable, conservative and irreducible single birth $Q$-matrix. The minimal process is recurrent if and only if $\sum_{\ell=0}^{+\infty}F_\ell^{(0)}=+\infty$.
\end{proposition}

To discuss uniqueness, define a $Q$-matrix $\overline{Q}=(\overline{q}_{ij}:i,j\in\overline{E})$ on the extended state space $\overline{E}=E\cup\{-1\}$ by
\[
\begin{aligned}
\overline{q}_{i,-1}&=1,\quad \overline{q}_{ii}=q_{ii}-1,\quad \overline{q}_{ij}=q_{ij},\quad j\neq i,\,i,j\in E;\\
\overline{q}_{-1,0}&=1,\quad \overline{q}_{-1,-1}=-1,\quad \overline{q}_{-1,j}=0,\quad 1\le j\in E.
\end{aligned}
\]
Then $\overline{Q}$ is irreducible and corresponds to a single birth process. Moreover, $\overline{T}_{-1,n+}=\sum_{0\le \ell\le n-1}\overline{F}_\ell^{(-1)}$, where $\overline{F}_\ell^{(-1)}$ admits an expression analogous to that of $F_\ell^{(n)}$ for the original matrix $Q$.
The following conclusion then follows from Theorem \ref{uni}.

\begin{theorem}
Let $Q=(q_{ij})$ be a totally stable and conservative single birth $Q$-matrix. The single birth process determined by $Q$ is unique if and only if $\sum_{\ell=-1}^{+\infty}\overline{F}_\ell^{(-1)}=+\infty$.
\end{theorem}

Note that
\[
\overline{q}_\ell^{(k)}=\sum_{j=-1}^{k}\overline{q}_{\ell j}=1+\sum_{j=0}^{k}q_{\ell j}=1+q_\ell^{(k)},\quad 0\le k<\ell;
\]
and
\[
\overline{F}_\ell^{(-1)}=\frac{1}{\overline{q}_{\ell,\ell+1}}\sum_{k=-1}^{\ell-1}\overline{q}_\ell^{(k)}\overline{F}_k^{(-1)}=\frac{1}{q_{\ell,\ell+1}}\left(\sum_{k=0}^{\ell-1}\big(1+q_\ell^{(k)}\big)\overline{F}_k^{(-1)}+1\right),\quad \ell\ge 0.
\]
Recall that $\widetilde{m}_\ell$ defined in \cite[Section 3]{C-Z-2014} satisfies
\[
\widetilde{m}_\ell=\frac{1}{q_{\ell,\ell+1}}\bigg(\sum_{k=0}^{\ell-1}\big(1+q_\ell^{(k)}\big)\widetilde{m}_k+1\bigg),\quad \ell\ge 0.
\]
Hence $\widetilde{m}_\ell=\overline{F}_\ell^{(-1)}$ for all $\ell\ge 0$. Meanwhile, \cite{C-Z-2014} also defines
\[
m_\ell=\frac{1}{q_{\ell,\ell+1}}\bigg(\sum_{0\le k\le \ell-1} q_\ell^{(k)}m_k+1\bigg),\quad \ell\ge 0,
\]
and proves that $\sum_{\ell=0}^{+\infty}\widetilde{m}_\ell=+\infty$ is equivalent to $\sum_{\ell=0}^{+\infty}m_\ell=+\infty$. Nevertheless, the equivalence between the recurrence of the minimal process on the extended space and the uniqueness of the original process was not recognized in \cite{C-Z-2014}. We may therefore restate the uniqueness criterion as follows (cf. \cite{Chen-1999, Y-C-1986}).

\begin{proposition}
Let $Q=(q_{ij})$ be a totally stable and conservative single birth $Q$-matrix. The single birth process determined by $Q$ is unique if and only if $\sum_{\ell=0}^{+\infty}m_\ell=+\infty$.
\end{proposition}
Therefore, even the above uniqueness criterion for birth-death processes is not obtained directly, but derived equivalently in this way. Of course, $m_\ell$ admits an intuitive probabilistic interpretation: $m_\ell = \mathbb E_\ell \sigma_{\ell+1}$. Moreover, by the single-birth property,
\[
R:=\sum_{\ell=0}^{+\infty}m_\ell=\lim_{n\to\infty}\sum_{\ell=0}^{n-1}m_\ell=\lim_{n\to\infty}\mathbb E_0\sigma_n =: \mathbb E_0\sigma_{\infty},
\]
and $R=\lim_{n\to\infty}\mathbb E_0\sigma_n=\lim_{n\to\infty}\mathbb E_0\sigma_{n+}=\mathbb E_0\eta$. Hence, the quantity $R$ possesses a nice probabilistic meaning: it is the mean time for the minimal process starting from state $0$ to first reach infinity, which is also the mean explosion time.

\section{Single Death Processes}
A process is called a single death process if $q_{i,i-1}>0$ for all $i\ge 1$, and $q_{i,i-j}=0$ for all $i\ge j\ge 2$. Define $q_i^{(k)}=\sum_{j=k}^{\infty} q_{ij}$ for $k>i\ge 0$, and set
\begin{equation}\label{D-Ex1}
G_n^{(n)}=1,\qquad G_i^{(n)}=\frac{1}{q_{i,i-1}}\sum_{k=i+1}^{n} q_i^{(k)} G_k^{(n)},\quad 0\le i<n.
\end{equation}
For notational convenience, we set \(q_{0,-1}=1\).

Parallel to Lemma \ref{C-Z-1}, we have the following lemma for single death processes.

\begin{lemma}(\cite[Corollary 2.3]{Z-2018})\label{Z-1}
Given $f$, the sequence $\{h_i\}$ defined recursively by
\[
h_i=\frac{1}{q_{i,i-1}}\bigg(\sum_{i+1\le k\le n} q_i^{(k)} h_k + f_i\bigg),\quad 0\le i\le n,
\]
admits the unified representation:
\[
h_i=\sum_{k=i}^{n} \frac{G_i^{(k)} f_k}{q_{k,k-1}},\quad 0\le i\le n.
\]
In particular, $\{G_n^{(k)}\}$ satisfies
\begin{equation}\label{D-Ex2}
G_n^{(n)}=1,\qquad G_i^{(n)}=\sum_{k=i}^{n-1} \frac{G_i^{(k)} q_k^{(n)}}{q_{k,k-1}},\quad 0\le i<n.
\end{equation}
\end{lemma}

\begin{proposition}\label{D-p1}
For a single death $Q$-matrix,  \eqref{eq0.1} has a unique solution
$$\aligned
&x_i=\sum_{i+1\le \ell\le n} G_\ell^{(n)} x_{n-1} - \sum_{i+1\le \ell\le n-1}\sum_{i\le k\le \ell-1} \frac{G_\ell^{(k)} V(k)}{q_{k,k-1}},\quad 0\le i\le n-2.\\
&x_{n-1}=\frac{1}{G_0^{(n)}}\sum_{0\le k\le n-1}\frac{G_0^{(k)}V(k)}{q_{k,k-1}}.
\endaligned$$
\end{proposition}

\begin{proof}
By \eqref{eq0.1} and the conservativeness of $Q$, for $0\le i\le n-1$, we obtain
\[\begin{aligned}
q_i^{(n)}x_i+\mathds 1_{\{1\le i\le n-1\}}q_{i,i-1}(x_i-x_{i-1})
&=\sum_{i+1\le j\le n-1} q_{ij}(x_j-x_i)+V(i) \\
%&=\sum_{i+1\le \ell\le n-1}(x_\ell-x_{\ell-1})\sum_{\ell\le j\le n-1} q_{ij}+V(i) \\
&=\sum_{i+1\le \ell\le n-1} q_i^{(\ell)}(x_\ell-x_{\ell-1})-q_i^{(n)}(x_{n-1}-x_i)+V(i).\\
\end{aligned}
\]
Therefore, for $1\le i\le n-1$, it follows from Lemma \ref{Z-1} and \eqref{D-Ex2} that
\[
\begin{aligned}
x_i-x_{i-1}
&=\frac{1}{q_{i,i-1}}\bigg(\sum_{i+1\le \ell\le n-1} q_i^{(\ell)}(x_\ell-x_{\ell-1})+V(i)-q_i^{(n)}x_{n-1}\bigg) \\
&=\sum_{i\le k\le n-1}\frac{G_i^{(k)}\big(V(k)-q_k^{(n)}x_{n-1}\big)}{q_{k,k-1}} \\
&=\sum_{i\le k\le n-1}\frac{G_i^{(k)}V(k)}{q_{k,k-1}}-G_i^{(n)}x_{n-1}.
\end{aligned}
\]
For $i=0$, we have
\[
0=\sum_{1\le \ell\le n-1} q_0^{(\ell)}(x_\ell-x_{\ell-1})-q_0^{(n)}x_{n-1}+V(0).
\]
Hence,
\[
\begin{aligned}
0&=\sum_{1\le \ell\le n-1} q_0^{(\ell)}\bigg(\sum_{\ell\le k\le n-1}\frac{G_\ell^{(k)}V(k)}{q_{k,k-1}}-G_\ell^{(n)}x_{n-1}\bigg)-q_0^{(n)}x_{n-1}+V(0) \\
&=\sum_{1\le k\le n-1}\frac{V(k)}{q_{k,k-1}}\sum_{1\le \ell\le k} q_0^{(\ell)}G_\ell^{(k)}-\sum_{1\le \ell\le n-1} q_0^{(\ell)}G_\ell^{(n)}x_{n-1}+V(0) \\
&=\sum_{0\le k\le n-1}\frac{G_0^{(k)}V(k)}{q_{k,k-1}}-G_0^{(n)}x_{n-1}
\end{aligned}
\]
by \eqref{D-Ex1}, i.e.,
\[
x_{n-1}=\frac{1}{G_0^{(n)}}\sum_{0\le k\le n-1}\frac{G_0^{(k)}V(k)}{q_{k,k-1}}.
\]
Therefore,  for $0\le i\le n-1$,
\[
\begin{aligned}
x_i &= x_{n-1}-\sum_{i+1\le \ell\le n-1}(x_\ell-x_{\ell-1}) \\
&= x_{n-1}-\sum_{i+1\le \ell\le n-1}\bigg(\sum_{\ell\le k\le n-1}\frac{G_\ell^{(k)}V(k)}{q_{k,k-1}}-G_\ell^{(n)}x_{n-1}\bigg) \\
&=\sum_{i+1\le \ell\le n}G_\ell^{(n)}x_{n-1}-\sum_{i+1\le \ell\le n-1}\sum_{\ell\le k\le n-1}\frac{G_\ell^{(k)}V(k)}{q_{k,k-1}}.
\end{aligned}
\]
This completes the proof.\end{proof}

Fix some $j\in\{0,1,\cdots, n-1\}$ and take $V(k)=\mathds 1_{\{k=j\}}$. Then by Proposition \ref{D-p1},
\[
x_{n-1}=\frac{G_0^{(j)}}{q_{j,j-1}G_0^{(n)}},\quad
x_i=
\begin{cases}
\displaystyle\sum_{i+1\le \ell\le n}G_\ell^{(n)}x_{n-1}, &0\le j\le i\le n-1;\\[6pt]
\displaystyle\sum_{i+1\le \ell\le n}G_\ell^{(n)}x_{n-1}-\sum_{i+1\le \ell<j}\frac{G_\ell^{(j)}}{q_{j,j-1}}, & 0\le i<j.
\end{cases}
\]
In particular, for $j=0$, we have $x_{n-1}=1/G_0^{(n)}$ and $\overline{T}_{0,n+}=x_0=\sum_{1\le \ell\le n}G_\ell^{(n)}/G_0^{(n)}$. Therefore, the following criterion follows directly from Theorem \ref{Rec} (cf. \cite{M-Y-Z-2026}).

\begin{proposition}
Let $Q=(q_{ij})$ be a totally stable, conservative and irreducible single death $Q$-matrix. The corresponding minimal process is recurrent if and only if
\[
\text{M-}\lim_{n\to\infty}\frac{\sum_{1\le \ell\le n}G_\ell^{(n)}}{G_0^{(n)}}=+\infty.
\]
\end{proposition}

To discuss the uniqueness of single death processes, following the discussion in Section 1, we define a new irreducible transition rate matrix on the extended state space $\overline{E}=\{-1,0,1,\dots\}$. For convenience of application, we consider a more general case directly: a totally stable, conservative and irreducible transition rate matrix on $\overline{E}$ given by
\begin{equation}\label{new-Q-2}
\begin{pmatrix}
-1 & \boldsymbol{ p}\\
\boldsymbol{ v}^\top & Q-\operatorname{diag}(\boldsymbol{ v})
\end{pmatrix}.
\end{equation}

Here $Q$ is a totally stable and conservative single death matrix on $E=\{0,1,2,\dots\}$, $\boldsymbol {v}=(v_0,v_1,v_2,\dots)$ is an infinite-dimensional nonnegative vector, $\operatorname{diag}(\boldsymbol{v})$ denotes the infinite matrix with $\boldsymbol{ v}$ on its main diagonal, $\boldsymbol{ p}=(p_0,p_1,p_2,\dots)$ is an infinite-dimensional nonnegative vector satisfying $\sum_{i\ge 0}p_i=1$. The choice of $\boldsymbol{ p}$ only needs to guarantee the irreducibility of the transition rate matrix. For simplicity, we still denote the entries of this matrix by $q_{ij}\;(i,j\in\overline{E})$.

Define
\[
\widetilde q_i^{(k)}=\sum_{j=k}^{\infty}q_{ij}+\mathds 1_{\{i\ge1\}}q_{i,-1}=q_i^{(k)}+\mathds 1_{\{i\ge1\}}q_{i,-1},\quad k>i\ge -1.
\]
Let
\[
\widetilde{G}_i^{(i)}=1,\qquad \widetilde{G}_i^{(n)}=\frac{1}{q_{i,i-1}}\sum_{k=i+1}^{n}\hat q_i^{(k)}\widetilde{G}_k^{(n)},\quad -1\le i<n.
\]
For notational convenience, we adopt the convention $q_{-1,-2}:=1$. We have \[
\sum_{i\le k\le n-1}\frac{\widetilde{G}_i^{(k)}\widetilde q_k^{(n)}}{q_{k,k-1}}=\widetilde{G}_i^{(n)}.\qquad\]Further set
\[
H_i^{(i)}=0,\qquad H_i^{(n)}=\frac{1}{q_{i,i-1}}\bigg(\sum_{i+1\le k\le n}\widetilde q_i^{(k)}H_k^{(n)}+\mathds{1}_{\{i\ge1\}}q_{i,-1}\bigg),\quad -1\le i\le n.
\]
Lemma \ref{Z-1} remains valid in the present setting. In particular,
\[
\sum_{i\le k\le n-1}\frac{\widetilde{G}_i^{(k)}q_{k,-1}\mathds{1}_{\{k\ge1\}}}{q_{k,k-1}}=H_i^{(n)},\quad -1\le i<n.
\]

For the above transition rate matrix on $\overline{E}$, we investigate the unique solution to the corresponding system \eqref{eq0.1}. Now \eqref{eq0.1} can be written as
\[
q_i x_i=\mathds{1}_{\{1\le i\le n-1\}}q_{i,-1}x_{-1}+\mathds{ 1}_{\{0\le i\le n-1\}}q_{i,i-1}x_{i-1}+\sum_{i+1\le j\le n-1}q_{ij}x_j+V(i),\quad -1\le i\le n-1.
\]
This, together with the conservativeness of $Q$, yields
\begin{equation}\label{D-Ex3}
\begin{aligned}
&q_i^{(n)}x_i+\mathds{1}_{\{1\le i\le n-1\}}q_{i,-1}x_i+\mathds{1}_{\{0\le i\le n-1\}}q_{i,i-1}(x_i-x_{i-1})\\
&=\mathds{1}_{\{1\le i\le n-1\}}q_{i,-1}x_{-1}+\sum_{i+1\le j\le n-1}q_{ij}(x_j-x_i)+V(i) \\
%&=\mathds{1}_{\{1\le i\le n-1\}}q_{i,-1}x_{-1}+\sum_{i+1\le \ell\le n-1}(x_\ell-x_{\ell-1})\sum_{\ell\le j\le n-1}q_{ij}+V(i) \\
&=\mathds{1}_{\{1\le i\le n-1\}}q_{i,-1}x_{-1}+\sum_{i+1\le \ell\le n-1} q_i^{(\ell)}(x_\ell-x_{\ell-1})-q_i^{(n)}(x_{n-1}-x_i)+V(i).
\end{aligned}
\end{equation}
Therefore, by Lemma \ref{Z-1}, for any $0\le i\le n-1$,
\[
\begin{aligned}
x_i-x_{i-1}
&=\frac{1}{q_{i,i-1}}\bigg(\sum_{i+1\le \ell\le n-1}\widetilde q_i^{(\ell)}(x_\ell-x_{\ell-1})+V(i)+\mathds{1}_{\{i\ge1\}}q_{i,-1}x_{-1}-\widetilde q_i^{(n)}x_{n-1}\bigg) \\
&=\sum_{i\le k\le n-1}\frac{\widetilde{G}_i^{(k)}\big(V(k)+q_{k,-1}\mathds{1}_{\{k\ge1\}}x_{-1}-\widetilde q_k^{(n)}x_{n-1}\big)}{q_{k,k-1}} \\
&=\sum_{i\le k\le n-1}\frac{\widetilde{G}_i^{(k)}V(k)}{q_{k,k-1}}+H_i^{(n)}x_{-1}-\widetilde{G}_i^{(n)}x_{n-1}.
\end{aligned}
\]
Summing the above identity from 0 to \(n-1\), we obtain
\[
x_{n-1}-x_{-1}=\sum_{0\le i\le n-1}\sum_{i\le k\le n-1}\frac{\widetilde{G}_i^{(k)}V(k)}{q_{k,k-1}}+\sum_{0\le i\le n-1}H_i^{(n)}x_{-1}-\sum_{0\le i\le n-1}\widetilde{G}_i^{(n)}x_{n-1},
\]
i.e.,
\begin{equation}\label{eq2.1}
\sum_{0\le i\le n}\widetilde{G}_i^{(n)}x_{n-1}-\bigg(\sum_{0\le i\le n}H_i^{(n)}+1\bigg)x_{-1}=\sum_{0\le i\le n-1}\sum_{i\le k\le n-1}\frac{\widetilde{G}_i^{(k)}V(k)}{q_{k,k-1}}. \end{equation}
Letting  $i=-1$ in the system \eqref{D-Ex3} , we obtain
\[
\begin{aligned}
0&=\sum_{0\le \ell\le n-1} q_{-1}^{(\ell)}(x_\ell-x_{\ell-1})+V(-1)- q_{-1}^{(n)}x_{n-1} \\
&=\sum_{0\le \ell\le n-1}\widetilde q_{-1}^{(\ell)}\bigg(\sum_{\ell\le k\le n-1}\frac{\widetilde{G}_\ell^{(k)}V(k)}{q_{k,k-1}}+H_\ell^{(n)}x_{\ell-1}-\widetilde{G}_\ell^{(n)}x_{n-1}\bigg)+V(-1)-\widetilde q_{-1}^{(n)}x_{n-1} \\
&=\sum_{-1\le k\le n-1}\frac{\widetilde{G}_{-1}^{(k)}V(k)}{q_{k,k-1}}+H_{-1}^{(n)}x_{-1}-\widetilde{G}_{-1}^{(n)}x_{n-1},
\end{aligned}
\]
i.e.,
\begin{equation}\label{eq2.2}
\widetilde{G}_{-1}^{(n)}x_{n-1}-H_{-1}^{(n)}x_{-1}=\sum_{-1\le k\le n-1}\frac{\widetilde{G}_{-1}^{(k)}V(k)}{q_{k,k-1}}.
\end{equation}
Solving the linear system consisting of equations \eqref{eq2.1} and \eqref{eq2.2}, we obtain
\[
x_{-1}=\frac{\widetilde{G}_{-1}^{(n)}\displaystyle\sum_{0\le i\le n-1}\sum_{i\le k\le n-1}\frac{\widetilde{G}_i^{(k)}V(k)}{q_{k,k-1}}-\sum_{0\le i\le n}\widetilde{G}_i^{(n)}\displaystyle\sum_{-1\le k\le n-1}\frac{\widetilde{G}_{-1}^{(k)}V(k)}{q_{k,k-1}}}{H_{-1}^{(n)}\displaystyle\sum_{0\le i\le n}\widetilde{G}_i^{(n)}-\widetilde{G}_{-1}^{(n)}\bigg(\sum_{0\le i\le n}H_i^{(n)}+1\bigg)}.
\]
Take $V(k)=\mathds{1}_{\{k=-1\}}$, then it follows directly that
\[
\overline{T}_{-1,n+}=x_{-1}=\frac{\displaystyle\sum_{0\le i\le n}\widetilde{G}_i^{(n)}}{\widetilde{G}_{-1}^{(n)}\bigg(\displaystyle\sum_{0\le i\le n}H_i^{(n)}+1\bigg)-H_{-1}^{(n)}\displaystyle\sum_{0\le i\le n}\widetilde{G}_i^{(n)}}.
\]
Hence, by Theorem \ref{Rec}, the minimal process on $\overline{E}$ corresponding to \eqref{new-Q-2} is recurrent if and only if $\overline{T}_{-1,n+}$ increases monotonically to $+\infty$ as $n\to\infty$.

Return to the uniqueness of single death processes. Let $\boldsymbol{v}=(1,1,1,\dots)$, $\widetilde q_i^{(k)}=q_i^{(k)}+\mathds{1}_{\{i\ge1\}}$ for $k>i\ge -1$, and
\[
H_i^{(n)}=\frac{1}{q_{i,i-1}}\bigg(\sum_{i+1\le k\le n}\widetilde q_i^{(k)}H_k^{(n)}+\mathds{1}_{\{i\ge1\}}\bigg),\quad -1\le i<n.
\]
Moreover, by Lemma \ref{Z-1},
\[
\sum_{i\le k\le n-1}\frac{\widetilde{G}_i^{(k)}\mathds{1}_{\{k\ge1\}}}{q_{k,k-1}}=H_i^{(n)},\quad -1\le i<n.
\]
We arrive at the following uniqueness criterion for single death processes.

\begin{theorem}
Let $Q=(q_{ij})$ be a totally stable and conservative single death $Q$-matrix. The single death process determined by $Q$ is unique if and only if
\[
\text{M-}\lim_{n\to\infty}\frac{\displaystyle\sum_{0\le i\le n}\widetilde{G}_i^{(n)}}{\widetilde{G}_{-1}^{(n)}\bigg(\displaystyle\sum_{0\le i\le n}H_i^{(n)}+1\bigg)-H_{-1}^{(n)}\displaystyle\sum_{0\le i\le n}\widetilde{G}_i^{(n)}}=+\infty.
\]
\end{theorem}

Furthermore, if we additionally assume that the single death matrix is irreducible, we may take $q_{-1}=q_{-1,0}=1$. Then $\widetilde{G}_{-1}^{(n)}=\widetilde{G}_0^{(n)}$ and $H_{-1}^{(n)}=H_0^{(n)}$. Therefore, the single death process is unique if and only if
\[
\widetilde{G}_0^{(n)}\cdot\frac{\displaystyle\sum_{0\le i\le n}H_i^{(n)}+1}{\displaystyle\sum_{0\le i\le n}\widetilde{G}_i^{(n)}}-H_0^{(n)}\downarrow 0,
\]
as $n\to\infty$.

In this setting, we compare our results with those in \cite{M-Y-Z-2026}, where $\overline{G}_k^{(n)}$ and $\overline{H}_k^{(n)}$ denote the quantities appearing in the uniqueness criterion.
The relations between these and our quantities are: for all $k\ge 1$, $\overline{G}_k^{(n)}=\widetilde{G}_k^{(n)}$ and $\overline{H}_k^{(n)}=H_k^{(n)}$, and $\overline{G}_0^{(n)}=\sum_{i=0}^n\widetilde{G}_i^{(n)}$ and $\overline{H}_0^{(n)}=\sum_{i=0}^nH_i^{(n)}+1$.
 Hence the above condition is equivalent to
\[
\text{M-}\lim_{n\to\infty}\bigg(\frac{\overline{H}_0^{(n)}}{\overline{G}_0^{(n)}}\sum_{1\le i\le n}\overline{G}_i^{(n)}-\sum_{1\le i\le n}\overline{H}_i^{(n)}\bigg)=1,
\]
which recovers the condition in \cite{M-Y-Z-2026}.

\section{Birth-death processes on trees}
A tree is an undirected, locally compact, connected acyclic graph. Let $T$ denote the set of all vertices of the tree, and $o\in T$ the unique root. For any non-root vertex $x$ on the tree, there is a unique simple path from $x$ to $o$. The set consisting of all non-root vertices on this path is denoted by $\mathscr{P}(x)$. The number of edges on this path is called the level of $x$, written $|x|$. Set $|o|=0$. Denote $E_k=\{x\in T:|x|=k\}$, which is called the $k$-th layer of the tree. Two vertices connected by an edge are said to be adjacent. A vertex adjacent to $x$ in the $(|x|+1)$-th layer is called a child of $x$. The set of all children of $x$ is denoted by $J(x)$. Assume that $J(x)$ is finite for every $x\in T$. The unique vertex adjacent to $x$ in the $(|x|-1)$-th layer is called the parent of $x$, denoted by $x^*$. $T_x$ stands for the subtree rooted at $x$ containing all descendants of $x$, as well as the vertex set of this subtree.

A birth-death process on a tree is a continuous-time symmetric random walk $X=\{X(t),t\ge 0\}$ on $T$, whose transition rate matrix $Q=(q(x,y))$ satisfies $q(x,y)>0$ if and only if $y=x^*$ or $y\in J(x)$. We adopt the convention that $q(o,o^*)=0$. Suppose that $Q$ is conservative and totally stable, i.e.,
\[
q(x):=-q(x,x)=q(x,x^*)+\sum_{y\in J(x)}q(x,y)<+\infty,\qquad x\in T.
\]

We now consider sufficiency criteria and uniqueness criteria for birth-death processes on trees. Note that in this tree-version setting, we define $T_{x,n+}(V)=T_{x,A}(V)$ for $0\le |x|<n$ and $\sigma_{n+}=\sigma_A$, where $A=\bigcup_{k=n}^{\infty}E_k$; we write $T_{x,n+}(V)$ as $T_{x,n+}$ when $V(x)=\mathds{1}_{\{x=o\}}$.
Moreover, $\{T_{x,n+}(V),\,0\le|x|<n\}$ is the minimal nonnegative solution (also the unique solution) to the following tree-version of system \eqref{eq0.1}:
\begin{equation}\label{eq5.1}
u(x)=\sum_{\substack{y\ne x\\|y|<n}}\frac{q(x,y)}{q(x)}u(y)+\frac{V(x)}{q(x)},\qquad 0\le|x|<n.
\end{equation}
 Define
\[
G_x^{(n)}=\sum_{y\in J(x)}q(x,y),\quad |x|=n-1;
\]
\[
G_x^{(n)}=\sum_{y\in J(x)}\frac{q(x,y)G_y^{(n)}}{q(y,x)+G_y^{(n)}},\qquad 0\le|x|<n-1.
\]

\begin{lemma}Let
\[
V_x^{(n)}=\frac{1}{q(x,x^*)+G_x^{(n)}}\bigg(\sum_{\substack{y\in J(x)\\|y|\le n}}q(x,y)V_y^{(n)}+V(x)\bigg),\qquad 1\le|x|\le n-1;
\]
\[
V_o^{(n)}=\frac{1}{G_o^{(n)}}\bigg(\sum_{y\in J(o)}q(o,y)V_y^{(n)}+V(o)\bigg).
\]
Then the unique solution to system \eqref{eq5.1} is given by
\begin{equation}\label{D-T-f1}
u(o)=V_o^{(n)},\quad
u(x)=\frac{q(x,x^*)}{q(x,x^*)+G_x^{(n)}}\cdot u(x^*)+V_x^{(n)},\qquad 1\le|x|<n.
\end{equation}\end{lemma}

\begin{proof}
System \eqref{eq5.1} admits the equivalent form
\[
\begin{cases}
q(o)u(o)=\displaystyle\sum_{y\in J(o)}q(o,y)u(y)+V(o),\\[6pt]
q(x)u(x)=q(x,x^*)u(x^*)+V(x), & |x|=n-1,\\[6pt]
q(x)u(x)=q(x,x^*)u(x^*)+\displaystyle\sum_{y\in J(x)}q(x,y)u(y)+V(x), & 1\le|x|<n-1.
\end{cases}
\]
For $|x|=n-1$, $q(x)=q(x,x^*)+G_x^{(n)}$, so
\[
u(x)=\frac{q(x,x^*)}{q(x,x^*)+G_x^{(n)}}\cdot u(x^*)+\frac{V(x)}{q(x,x^*)+G_x^{(n)}}
=\frac{q(x,x^*)}{q(x,x^*)+G_x^{(n)}}\cdot u(x^*)+V_x^{(n)}.
\]
Thus \eqref{D-T-f1} holds for $|x|=n-1$. Assume that \eqref{D-T-f1}  is true for all vertices satisfying $k+1\le|x|<n$. Consider $|x|=k\ge 1$. Then
\[
\begin{aligned}
q(x)u(x)
&=q(x,x^*)u(x^*)+\sum_{y\in J(x)}q(x,y)u(y)+V(x)\\
&=q(x,x^*)u(x^*)+\sum_{y\in J(x)}q(x,y)\bigg(\frac{q(y,y^*)}{q(y,y^*)+G_y^{(n)}}\cdot u(y^*)+V_y^{(n)}\bigg)+V(x)\\
&=q(x,x^*)u(x^*)+\sum_{y\in J(x)}q(x,y)\bigg(\frac{q(y,y^*)}{q(y,y^*)+G_y^{(n)}}\cdot u(x)+V_y^{(n)}\bigg)+V(x).
\end{aligned}
\]
By the conservativeness of $Q$, rearranging the above identity yields
\[
\bigg(q(x,x^*)+\sum_{y\in J(x)}\frac{q(x,y)G_y^{(n)}}{q(y,y^*)+G_y^{(n)}}\bigg)u(x)
=q(x,x^*)u(x^*)+\sum_{y\in J(x)}q(x,y)V_y^{(n)}+V(x).
\]
That is,
\[
\big(q(x,x^*)+G_x^{(n)}\big)u(x)=q(x,x^*)u(x^*)+\sum_{y\in J(x)}q(x,y)V_y^{(n)}+V(x).
\]
Therefore,
\[
\begin{aligned}
u(x)&=\frac{q(x,x^*)}{q(x,x^*)+G_x^{(n)}}\cdot u(x^*)+\frac{1}{q(x,x^*)+G_x^{(n)}}\bigg(\sum_{y\in J(x)}q(x,y)V_y^{(n)}+V(x)\bigg)\\
&=\frac{q(x,x^*)}{q(x,x^*)+G_x^{(n)}}\cdot u(x^*)+V_x^{(n)}.
\end{aligned}
\]
Thus \eqref{D-T-f1} holds for $|x|=k\ge 1$. By induction, \eqref{D-T-f1}  holds for all $1\le |x|<n$.
For $x=o$,
\[
\begin{aligned}
q(o)u(o)
&=\sum_{y\in J(o)}q(o,y)\bigg(\frac{q(y,y^*)}{q(y,y^*)+G_y^{(n)}}\cdot u(y^*)+V_y^{(n)}\bigg)+V(o)\\
&=\sum_{y\in J(o)}q(o,y)\bigg(\frac{q(y,y^*)}{q(y,y^*)+G_y^{(n)}}\cdot u(o)+V_y^{(n)}\bigg)+V(o).
\end{aligned}
\]
Rearranging and using $q(o)=\sum_{y\in J(o)}q(o,y)$, we obtain
\[
\sum_{y\in J(o)}\frac{q(o,y)G_y^{(n)}}{q(y,y^*)+G_y^{(n)}}\cdot u(o)
=\sum_{y\in J(o)}q(o,y)V_y^{(n)}+V(o),
\]
i.e.,
\[
G_o^{(n)}u(o)=\sum_{y\in J(o)}q(o,y)V_y^{(n)}+V(o).
\]
Hence $u(o)=V_o^{(n)}$. This completes the proof.
\end{proof}

Next, take $V(x)=\mathds{1}_{\{x=o\}}$. Then $V_x^{(n)}=0$ for $x\neq o$, and $V_o^{(n)}=1/G_o^{(n)}$. Therefore $T_{o,n+}=u(o)=1/G_o^{(n)}$. Then by Theorem \ref{Rec} we directly obtain the following criterion (see \cite{W-Y-Z-2025}).

\begin{proposition}
Let $Q=(q(x,y))$ be a totally stable and conservative birth-death $Q$-matrix on a tree. The corresponding minimal process is recurrent if and only if
\[
\text{M-}\lim_{n\to\infty}\frac{1}{G_o^{(n)}}=+\infty,
\]
equivalently, $\text{M-}\lim_{n\to\infty}G_o^{(n)}=0$.
\end{proposition}

To study the uniqueness of birth-death processes on trees, we define a new irreducible transition rate matrix $\overline{Q}=(\bar q(x,y))$ on the enlarged state space $T\cup\{-1\}$ as follows:
\[
\bar q(-1,o)=1,\qquad \bar q(-1,x)=0,\quad o\neq x\in T,
\]
and
\[
\bar q(x,-1)=1,\qquad \bar q(x,y)=q(x,y),\quad x\neq y,\;x,y\in T.
\]
To keep conservativeness, set
\[
\bar q(-1)=-\bar q(-1,-1)=1,\quad \bar q(x)=-\bar q(x,x)=1-q(x,x)=1+q(x),\quad x\in T.
\]

For the transition rate matrix $\overline{Q}=(\bar q(x,y))$ on $T\cup\{-1\}$, we consider the unique nonnegative solution to tree-version of system \eqref{eq0.1}.
\iffalse
\begin{cases}
\bar q(-1)u(-1)=\bar q(-1,o)u(o)+V(-1),\\[4pt]
\bar q(o)u(o)=\bar q(o,-1)u(-1)+\displaystyle\sum_{y\in J(o)}\bar q(o,y)u(y)+V(o),\\[6pt]
\bar q(x)u(x)=\bar q(x,-1)u(-1)+\bar q(x,x^*)u(x^*)+\displaystyle\sum_{y\in J(x)}\bar q(x,y)u(y)+V(x), & 1\le|x|<n-1,\\[6pt]
\bar q(x)u(x)=\bar q(x,-1)u(-1)+\bar q(x,x^*)u(x^*)+V(x), & |x|=n-1.
\end{cases}
\]
i.e.,
\[
\begin{cases}
u(-1)=u(o)+V(-1),\\[4pt]
\big(1+q(o)\big)u(o)=u(-1)+\displaystyle\sum_{y\in J(o)}q(o,y)u(y)+V(o),\\[6pt]
\big(1+q(x)\big)u(x)=u(-1)+q(x,x^*)u(x^*)+\displaystyle\sum_{y\in J(x)}q(x,y)u(y)+V(x), & 1\le|x|<n-1,\\[6pt]
\big(1+q(x)\big)u(x)=u(-1)+q(x,x^*)u(x^*)+V(x), & |x|=n-1.
\end{cases}
\]
\fi
Let $\overline{u}=u-u(-1)$. Then the system \eqref{eq0.1} can be wirtten as
\begin{equation}\label{eq5.2}
\begin{cases}
0=\overline{u}(o)+V(-1),\\[4pt]
\big(1+q(o)\big)\overline{u}(o)=\displaystyle\sum_{y\in J(o)}q(o,y)\overline{u}(y)+V(o),\\[6pt]
\big(1+q(x)\big)\overline{u}(x)=q(x,x^*)\overline{u}(x^*)+\displaystyle\sum_{y\in J(x)}q(x,y)\overline{u}(y)+V(x), & 1\le|x|<n-1,\\[6pt]
\big(1+q(x)\big)\overline{u}(x)=q(x,x^*)\overline{u}(x^*)-\displaystyle\sum_{y\in J(x)}q(x,y)u(-1)+V(x), & |x|=n-1.
\end{cases}
\end{equation}

Define
\[
\overline{G}_x^{(n)}=\sum_{y\in J(x)}q(x,y)+1,\quad |x|=n-1;\qquad
\overline{G}_x^{(n)}=\sum_{y\in J(x)}\frac{q(x,y)\overline{G}_y^{(n)}}{q(y,x)+\overline{G}_y^{(n)}}+1,\quad 0\le|x|<n-1;
\]
\[
H_x^{(n)}=\sum_{y\in J(x)}q(x,y),\quad |x|=n-1;\qquad
H_x^{(n)}=\sum_{y\in J(x)}\frac{q(x,y)H_y^{(n)}}{q(y,x)+\overline{G}_y^{(n)}},\quad 0\le|x|<n-1.
\]
\begin{lemma}Let
\[
\overline{V}_x^{(n)}=\frac{1}{q(x,x^*)+\overline{G}_x^{(n)}}\bigg(\sum_{\substack{y\in J(x)\\|y|\le n}}q(x,y)\overline{V}_y^{(n)}+V(x)\bigg),\qquad 1\le|x|\le n-1;
\]
\[
\overline{V}_o^{(n)}=\frac{1}{\overline{G}_o^{(n)}}\bigg(\sum_{y\in J(o)}q(o,y)\overline{V}_y^{(n)}+V(o)\bigg).
\]
Then the unique solution to system \eqref{eq5.2} satisfies
\[
u(-1)=\frac{\overline{G}_o^{(n)}}{H_o^{(n)}}\big(\overline{V}_o^{(n)}+V(-1)\big),\quad
\overline{u}(o)=-\frac{H_o^{(n)}}{\overline{G}_o^{(n)}}\cdot u(-1)+\overline{V}_o^{(n)},
\]
\begin{equation}\label{D-T-f2}
\overline{u}(x)=\frac{q(x,x^*)}{q(x,x^*)+\overline{G}_x^{(n)}}\cdot\overline{u}(x^*)
-\frac{H_x^{(n)}}{q(x,x^*)+\overline{G}_x^{(n)}}+\overline{V}_x^{(n)},\qquad 1\le|x|<n.
\end{equation}
\end{lemma}

\begin{proof}
For $|x|=n-1$, we have $1+q(x)=q(x,x^*)+\overline{G}_x^{(n)}$. Hence
\[
\begin{aligned}
\overline{u}(x)
&=\frac{q(x,x^*)}{q(x,x^*)+\overline{G}_x^{(n)}}\cdot\overline{u}(x^*)
-\frac{\displaystyle\sum_{y\in J(x)}q(x,y)u(-1)}{q(x,x^*)+\overline{G}_x^{(n)}}
+\frac{V(x)}{q(x,x^*)+\overline{G}_x^{(n)}}\\
&=\frac{q(x,x^*)}{q(x,x^*)+\overline{G}_x^{(n)}}\cdot\overline{u}(x^*)
-\frac{H_x^{(n)}u(-1)}{q(x,x^*)+\overline{G}_x^{(n)}}+\overline{V}_x^{(n)}.
\end{aligned}
\]
Thus \eqref{D-T-f2} holds for $|x|=n-1$. Assume that \eqref{D-T-f2} is true for all vertices satisfying $k+1\le|x|<n$. Consider $|x|=k\ge 1$. Then
\[
\begin{aligned}
&\big(1+q(x)\big)\overline{u}(x)
=q(x,x^*)\overline{u}(x^*)+\sum_{y\in J(x)}q(x,y)\overline{u}(y)+V(x)\\
&=q(x,x^*)\overline{u}(x^*)+\sum_{y\in J(x)}q(x,y)\bigg(\frac{q(y,y^*)}{q(y,y^*)+\overline{G}_y^{(n)}}\overline{u}(y^*)
-\frac{H_y^{(n)}u(-1)}{q(y,y^*)+\overline{G}_y^{(n)}}+\overline{V}_y^{(n)}\bigg)+V(x)\\
&=q(x,x^*)\overline{u}(x^*)+\sum_{y\in J(x)}q(x,y)\bigg(\frac{q(y,y^*)}{q(y,y^*)+\overline{G}_y^{(n)}}\overline{u}(x)
+\overline{V}_y^{(n)}\bigg)-H_x^{(n)}u(-1)+V(x).
\end{aligned}
\]
Rearranging terms, and using $q(x)=q(x,x^*)+\sum_{y\in J(x)}q(x,y)$, $1\le |x|<n$, we obtain
\[
\big(q(x,x^*)+\overline{G}_x^{(n)}\big)\overline{u}(x)
=q(x,x^*)\overline{u}(x^*)-H_x^{(n)}u(-1)+\sum_{y\in J(x)}q(x,y)\overline{V}_y^{(n)}+V(x).
\]
Therefore,
\[
\overline{u}(x)=\frac{q(x,x^*)}{q(x,x^*)+\overline{G}_x^{(n)}}\cdot\overline{u}(x^*)
-\frac{H_x^{(n)}u(-1)}{q(x,x^*)+\overline{G}_x^{(n)}}+\overline{V}_x^{(n)}.
\]
Hence \eqref{D-T-f2} holds for $|x|=k\ge 1$. By induction, \eqref{D-T-f2} holds for all $1\le|x|<n$. Moreover,
\[
\begin{aligned}
\big(1+q(o)\big)\overline{u}(o)
&=\sum_{y\in J(o)}q(o,y)\bigg(\frac{q(y,y^*)}{q(y,y^*)+\overline{G}_y^{(n)}}\overline{u}(y^*)
-\frac{H_y^{(n)}u(-1)}{q(y,y^*)+\overline{G}_y^{(n)}}+\overline{V}_y^{(n)}\bigg)+V(o)\\
&=\sum_{y\in J(o)}q(o,y)\bigg(\frac{q(y,y^*)}{q(y,y^*)+\overline{G}_y^{(n)}}\overline{u}(o)
-\frac{H_y^{(n)}u(-1)}{q(y,y^*)+\overline{G}_y^{(n)}}+\overline{V}_y^{(n)}\bigg)+V(o).
\end{aligned}
\]
Rearranging and using $q(o)=\sum_{y\in J(o)}q(o,y)$, we obtain
\[
\overline{G}_o^{(n)}\overline{u}(o)=-H_o^{(n)}u(-1)+\sum_{y\in J(o)}q(o,y)\overline{V}_y^{(n)}+V(o).
\]
Thus $\overline{u}(o)=-H_o^{(n)}u(-1)/\overline{G}_o^{(n)}+\overline{V}_o^{(n)}$. Meanwhile, from system \eqref{eq5.2}, $\overline{u}(o)=-V(-1)$. It follows that
\[
u(-1)=\overline{G}_o^{(n)}\big(\overline{V}_o^{(n)}+V(-1)\big)/H_o^{(n)}.
\]
The proof of the lemma is complete.
\end{proof}

Next, take $V(x)=\mathds{1}_{\{x=-1\}}$. Then $\overline{V}_x^{(n)}=0$ for $x\in T$. Hence $\overline{T}_{-1,n+}=u(-1)=\overline{G}_o^{(n)}/H_o^{(n)}$. Therefore, using Theorem \ref{uni}, we directly recover the following criterion from \cite{W-Y-Z-2025}.
\begin{proposition}
Let $Q=(q(x,y))$ be a totally stable and conservative birth-death $Q$-matrix on a tree. The process determined by $Q$ is unique if and only if
\[
\text{M-}\lim_{n\to\infty}\frac{\overline{G}_o^{(n)}}{H_o^{(n)}}=+\infty,
\]
or equivalently,
\[
\text{M-}\lim_{n\to\infty}\frac{H_o^{(n)}}{\overline{G}_o^{(n)}}=0.\]
\end{proposition}

\end{document}